\documentclass{amsart}
\usepackage{graphicx,amsfonts,amsmath,amssymb,amsthm,dsfont}
\newtheorem{thm}{Theorem}[section]

\newtheorem{lem}[thm]{Lemma}

\theoremstyle{definition}

\theoremstyle{remark}
\newtheorem{rem}[thm]{Remark}
\numberwithin{equation}{section}

\newcommand{\eps}{\varepsilon}

\newcommand{\A}{\mathcal{A}}
\newcommand{\cD}{\mathcal{D}}
\newcommand{\tsq}{\dbnframe\hspace{-8.25pt}\smoothing}

\newcommand{\Q}{{\mathbb Q}}
\newcommand{\tetrahedron}{\begin{array}{c}\begin{picture}(15,20)
\put(-3,-2){\includegraphics[scale=0.26]{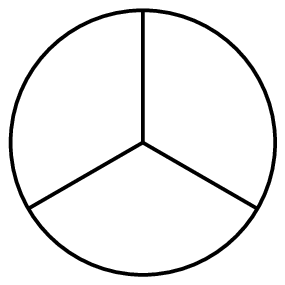}}\end{picture}\end{array}}
\renewcommand{\tsq}{\begin{array}{c}\begin{picture}(15,20)
\put(-3,-2){\includegraphics[scale=0.26]{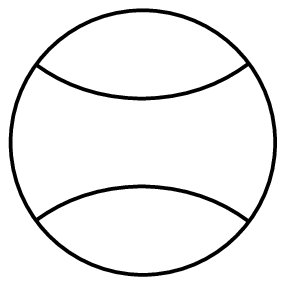}}\end{picture}\end{array}}
\newcommand{\wtr}{\begin{array}{c}\begin{picture}(20,20)
\put(-3,-2){\includegraphics[scale=0.26]{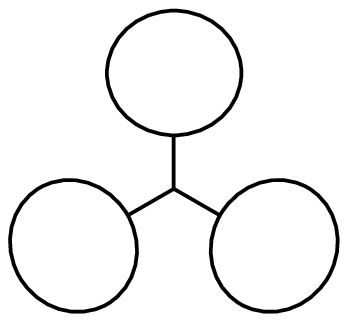}}\end{picture}\end{array}}
\newcommand{\mdl}{\begin{array}{c}\begin{picture}(24,20)
\put(-3,0){\includegraphics[scale=0.26]{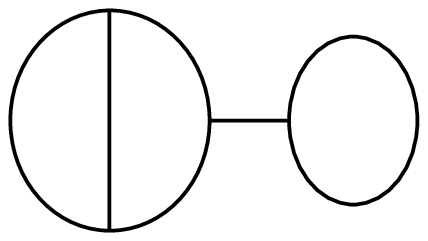}}\end{picture}\end{array}}
\newcommand{\bbl}{\begin{array}{c}\begin{picture}(24,20)
\put(-3,2){\includegraphics[scale=0.26]{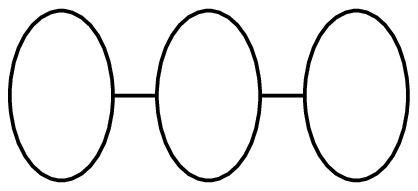}}\end{picture}\end{array}}
\newcommand{\pc}[2]{\mbox{$\begin{array}{c}
   \includegraphics[scale=#2]{#1.ps}
   \end{array}$}}
\newcommand{\vspc}[1]{\begin{picture}(0,#1)\end{picture}}

\allowdisplaybreaks

\begin{document}

\title[Vanishing of $3$--loop Jacobi diagrams of odd degree]
{Vanishing of $3$--loop Jacobi diagrams of odd degree}%

\author{Daniel Moskovich}%
\address{Research Institute for Mathematical Sciences,
Kyoto University, Kyoto, 606-8502 JAPAN}%
\email{dmoskovich@gmail.com}%
\urladdr{http://www.sumamathematica.com/}

\author{Tomotada Ohtsuki}%
\address{Research Institute for Mathematical Sciences,
Kyoto University, Kyoto, 606-8502 JAPAN}%
\email{tomotada@kurims.kyoto-u.ac.jp}%
\urladdr{http://www.kurims.kyoto-u.ac.jp/\~{}tomotada/}


\date{8th of August, 2006}%
\begin{abstract}
We prove the vanishing of the space of $3$--loop Jacobi diagrams of
odd degree. This implies that no $3$--loop Vassiliev invariant can
distinguish between a knot and its inverse.
\end{abstract}
\maketitle
\section{Introduction}

A Jacobi diagram is a uni-trivalent graph with some extra structure.
Such diagrams play a leading role in the theory of Vassiliev
invariants and Kontsevich invariants of knots. Vassiliev invariants
are defined by a filtration of the vector space spanned by knots,
whose graded spaces are identified with vector spaces spanned by
Jacobi diagrams subject to certain defining relations. The
Kontsevich invariant of a knot is defined as an infinite linear sum
of Jacobi diagrams. The physical background of these invariants is
in the perturbative expansion of the Chern--Simons path integral,
which is formulated in terms of uni-trivalent graphs; this is one
explanation why Jacobi diagrams appear in this theory. The
Kontsevich invariant is expected to classify knots, and from this
point of view it is important to identify the vector space spanned
by Jacobi diagrams subject to the defining relations.

It is conjectured that the space of Jacobi diagrams with an odd
number of legs vanishes; \cite{BN95b,problem_list}. This would imply
the claim that no Vassiliev invariant can distinguish a knot from
its inverse, where the {\it inverse} of an oriented knot is the knot
with the opposite orientation. In general, a knot and its inverse
are not isotopic, the simplest counter-example being the knot
$8_{17}$ with its two possible orientations. The consequences of the
possibility that Vassiliev invariants cannot make this distinction
are discussed in \cite{Kup96}. For the Lie algebra version of this
claim, see Remark \ref{rem.q_inv}. Dasbach claimed to have proved
the vanishing of $n$--loop Jacobi diagrams with an odd number of
legs for $n\leq 6$, but his proof has a gap for $n\geq 3$; see
Remark \ref{rem.Dasbach_gap}.

In the present paper, we prove the vanishing of $3$--loop Jacobi
diagrams with an odd number of legs (Theorem \ref{thm.main}). In our
proof, we consider the internal graph of a Jacobi diagram, which is
the trivalent graph obtained from the Jacobi diagram by removing its
legs, where a leg of a Jacobi diagram is an edge adjacent to a
univalent vertex. Then, following Nakatsuru \cite{Nts98}, we
identify each Jacobi diagram with a polynomial whose variables
correspond to the edges of the internal graph of the Jacobi diagram,
and present the space of $3$--loop Jacobi diagrams as a quotient
space of a direct sum of polynomial algebras corresponding to
$3$--loop internal graphs. Here, the quotient is derived from the
defining relations of Jacobi diagrams and from the symmetries of the
internal graphs. Thus, the proof is reduced to calculating the image
of the relations by the (skew) symmetrizer corresponding to the
internal graph's symmetry. This approach provides in passing
an alternative proof of \cite[Theorem 7.4]{Das98} in the `even
number of legs' case as well. The $4$--loop, $5$--loop, and
$6$--loop cases which Dasbach's result would have covered remain
open. In these higher loop degrees, the techniques used here lead to
more complicated calculations, which we have not been able to
complete. New ideas seem necessary in order to make further
progress.

The paper is organized as follows. In Section \ref{sec.Jd}, we
review several definitions concerning Jacobi diagrams and related
notions. In Section \ref{sec.3l_J_d}, we show how to identify the
space of $3$--loop Jacobi diagrams with a quotient space of a direct
sum of polynomial algebras and prove the vanishing of $3$--loop
Jacobi diagrams with an odd number of legs, which is the main
theorem of this paper. This proof requires the use of a certain
lemma, which we prove in Section \ref{sec.pf_lem}.

The gap in the proof of \cite[Theorem 5.4.3(\textrm{iii})]{Das97}
was discovered in a seminar when the first author tried to
generalize Dasbach's proof. The authors thank the participants of
the seminar --- Kazuo Habiro, Tadayuki Watanabe, and Atsushi Ishii
for their attention. The authors would especially like to thank
Pierre Vogel for useful comments regarding the identification of the
space of $n$--loop Jacobi diagrams. The first author would also like
to thank Alexander Stoimenow for useful discussions regarding
Dasbach's papers, and Oliver Dasbach for useful discussions.
The authors would also like to thank
the referees for their careful comments.

\section{Jacobi diagrams}
\label{sec.Jd}

In this section we review definitions of Jacobi diagrams, the space
of Jacobi diagrams, $n$--loop Jacobi diagrams, and define some
notations. For general references on the theory of Jacobi diagrams
see {\it e.g.} \cite{BN95b,Oht02}.

A {\it Jacobi diagram} is a graph whose vertices have valence $1$ or
$3$ and whose trivalent vertices are oriented {\it i.e.}, a cyclic
order of $3$ edges around each trivalent vertex is fixed. The {\it
degree} of a Jacobi diagram is defined to be half the total number
of vertices of the diagram. The {\it space of Jacobi diagrams} is
the vector space over $\Q$ spanned by Jacobi diagrams subject to the
AS (Anti--Symmetry) and IHX (written as ``I''$=$``H''$-$``X'')
relations, which are local moves between Jacobi diagrams which
differ inside a dotted circle as indicated below. The space of
Jacobi diagrams is graded by degree. (A Jacobi diagram of the type
we have just defined is sometimes called an {\it open Jacobi
diagram}, and the space of these Jacobi diagrams is sometimes
denoted $\mathcal{B}$ in the literature.)
\begin{enumerate}
\item[The {\it AS relation} \hspace*{-3.8pc}]$$
\begin{minipage}{42pt}
    \includegraphics[width=42pt]{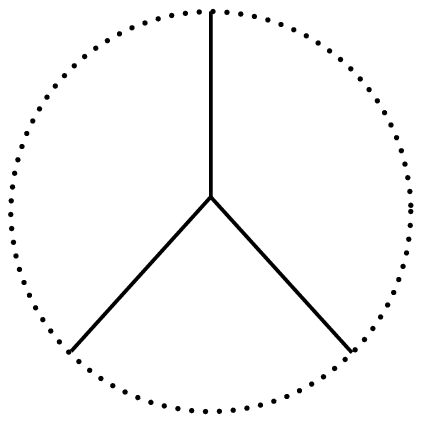}
\end{minipage}
\hspace{10pt}=\hspace{10pt}-\hspace{5pt}
\begin{minipage}{42pt}
    \includegraphics[width=42pt]{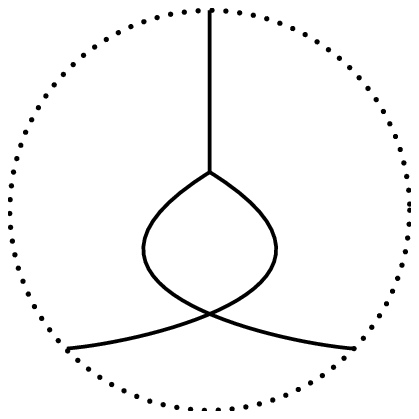}
\end{minipage}
$$
\item[The {\it IHX relation} \hspace*{-4.3pc}]$$
\begin{minipage}{42pt}
    \includegraphics[width=42pt]{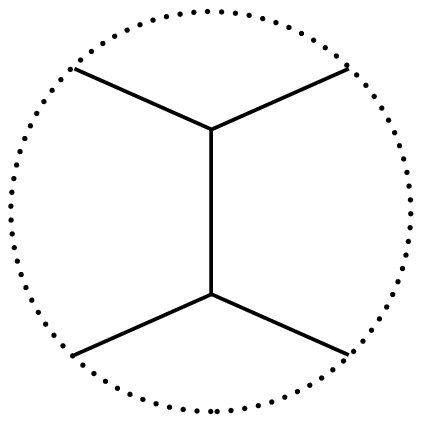}
\end{minipage}
\hspace{10pt}=\hspace{10pt}
\begin{minipage}{42pt}
    \includegraphics[width=42pt]{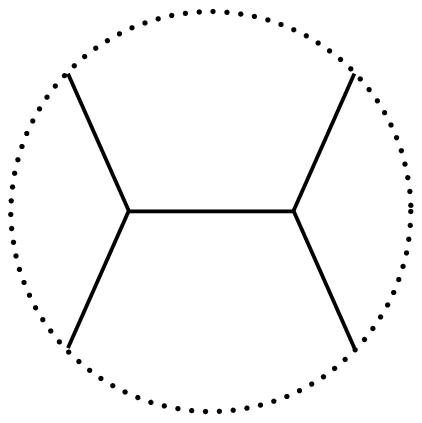}
\end{minipage}
\hspace{10pt}-\hspace{10pt}
\begin{minipage}{42pt}
    \includegraphics[width=42pt]{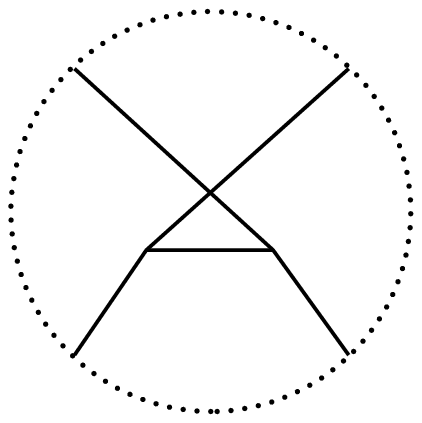}
\end{minipage}
$$
\end{enumerate}

A Jacobi diagram is called {\it $n$--loop} if
it is connected and its Euler number is
equal to $1-n$; {\it i.e.}, its first Betti number is equal to $n$.
(An $n$--loop Jacobi diagram is sometimes said to be of \emph{loop
degree $n-1$} in the literature.) We denote by $\A_{\,
\mbox{\scriptsize $n$--loop}}$ the space of $n$--loop Jacobi
diagrams, {\it i.e.}, the vector space spanned by $n$--loop Jacobi
diagrams subject to the AS and IHX relations. An edge adjacent to a
univalent vertex is called a {\it leg}.
We assume without loss of generality that
a Jacobi diagram does not have a trivalent vertex
which is adjacent to 2 legs,
since a Jacobi diagram with such a trivalent vertex vanishes
by the AS relation.
The {\it internal graph} of a Jacobi diagram is
the trivalent graph obtained from the Jacobi diagram by removing its legs.
We denote by $\A(\Gamma)$ the space of Jacobi diagrams
whose internal graph is $\Gamma$
modulo the action of the symmetry of $\Gamma$.

\section{$3$--loop Jacobi diagrams}
\label{sec.3l_J_d}

In this section we identify the space of $3$--loop Jacobi diagrams
as a graded vector space. In Section \ref{sec.s3lJd} we present the
space of $3$--loop Jacobi diagrams in terms of spaces $\A(\Gamma)$
for $3$--loop trivalent graphs $\Gamma$. In Section
\ref{sec.polyn_Jd} we present the space of such diagrams using
polynomial algebras. Using this presentation, we prove in Section
\ref{sec.odd_deg} that the odd degree part of this space vanishes,
which is the main theorem of this paper. In Section
\ref{sec.even_deg} we identify the even part of this space with some
polynomial algebra (following \cite{Nts98}).

\subsection{The space of $3$--loop Jacobi diagrams}
\label{sec.s3lJd}

In this section, we present the space of $3$--loop Jacobi diagrams
in terms of spaces $\A(\Gamma)$ for $3$--loop trivalent graphs
$\Gamma$.

Ignoring orientations of internal vertices, the internal graph of a
$3$--loop Jacobi diagram may be one of the five graphs below,
\begin{equation}
\label{eq.5ig}
\pc{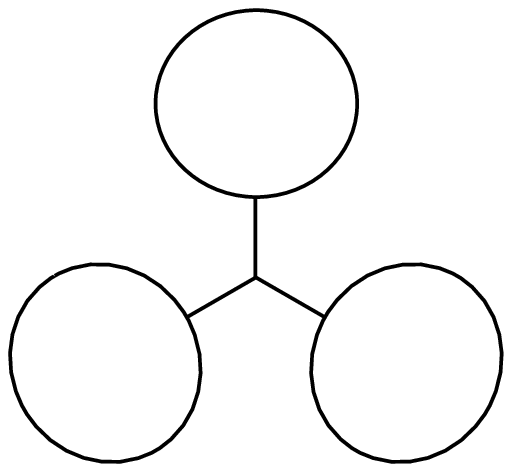}{0.33}, \
\pc{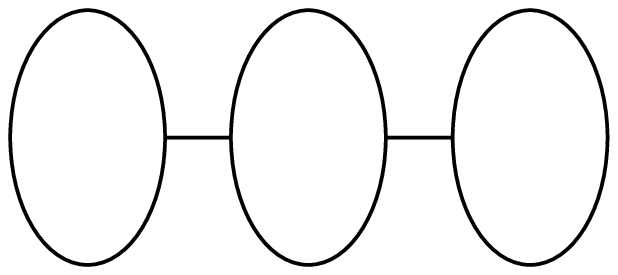}{0.33}, \
\pc{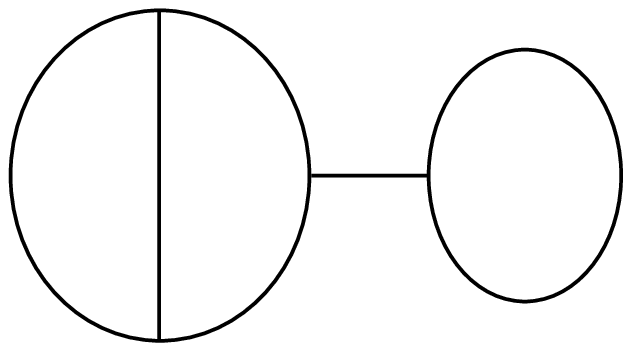}{0.33}, \
\pc{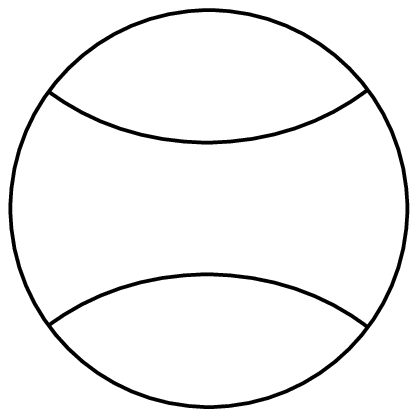}{0.33}, \
\pc{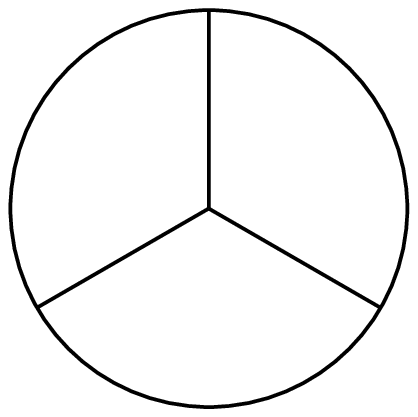}{0.33}.
\end{equation}

The space of $3$--loop Jacobi diagrams is presented by
\begin{equation}
\label{eq.A3l}
\A_{\, \mbox{\scriptsize 3--loop}}  \ \cong \
\Big( \bigoplus_{\mbox{\scriptsize $\Gamma$ in (\ref{eq.5ig})}} \!\!
\A (\Gamma) \Big) \Big/ \, {\rm IHX} ,
\end{equation}
where ``IHX'' implies the IHX relations among these $\Gamma$; all
such relations are obtained by replacing a neighborhood of a
$4$--valent vertex of one of the following graphs with the defining
graphs of the IHX relation,
\begin{equation}
\label{eq.4vg_ihx}
\pc{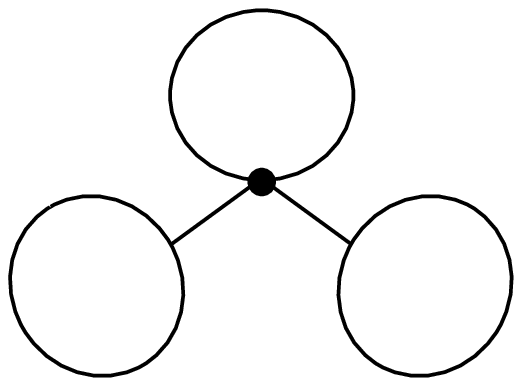}{0.32}, \
\pc{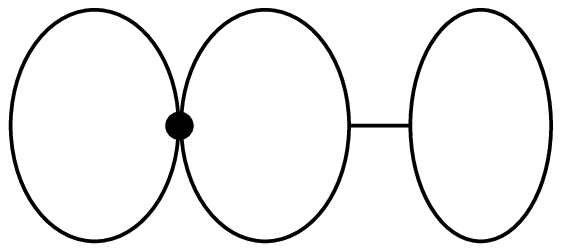}{0.32}, \
\pc{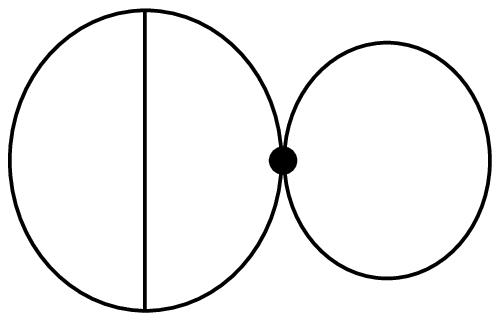}{0.32}, \
\pc{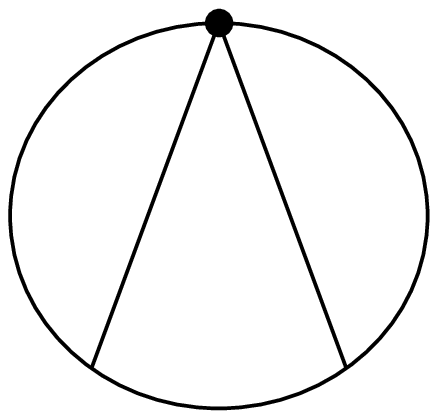}{0.32}, \
\pc{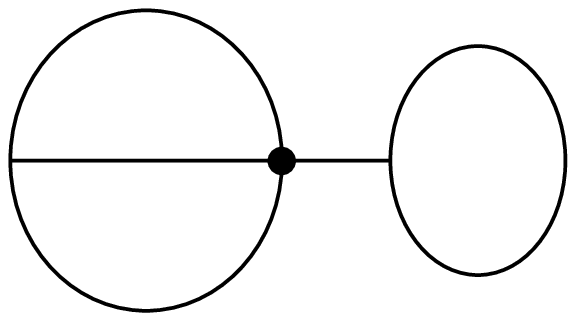}{0.32}.
\end{equation}

We will see,
in Sections \ref{sec.odd_deg} and \ref{sec.even_deg}
for the odd and even degree parts respectively,
that (\ref{eq.A3l}) is isomorphic to
\begin{equation}
\label{eq.tsq_tet}
\A_{\, \mbox{\scriptsize 3--loop}}  \ \cong \ \Big( \A \big( \tsq
\big) \oplus \A \big( \tetrahedron \big) \Big) \Big/ \,
{\rm IHX},
\end{equation}
where this ``IHX'' implies the IHX relation obtained from the fourth
graph of (\ref{eq.4vg_ihx}). We describe $\A \big( \tsq \big)$ and
$\A \big( \tetrahedron \big)$ in terms of polynomial algebras in the
next section.

\subsection{Polynomial presentation of $3$--loop Jacobi diagrams}
\label{sec.polyn_Jd}

In this section we see that the space of $3$--loop Jacobi diagrams
is identified, as a graded vector space, with a quotient space of a
direct sum of polynomial algebras.

We identify $\A \big(\tetrahedron \big)$ with the polynomial algebra
on six letters signifying legs on each of the arcs of the internal
graphs, modulo the IHX relations on the legs, and modulo the action
of $\mathfrak{S}_{4}$ the automorphism group of the tetrahedron.
Thus:
$$
\A(\tetrahedron) \ \cong \ \Q[x_1,x_2,x_3,x_4,x_5,x_6] \big/
(\ref{eq.Atet_rel}), {\mathfrak S}_4 ,
$$
where
\begin{equation}
\label{eq.tet_polyn}
\begin{array}{c}
\begin{picture}(100,90)
\put(0,0){\includegraphics[width=90pt]{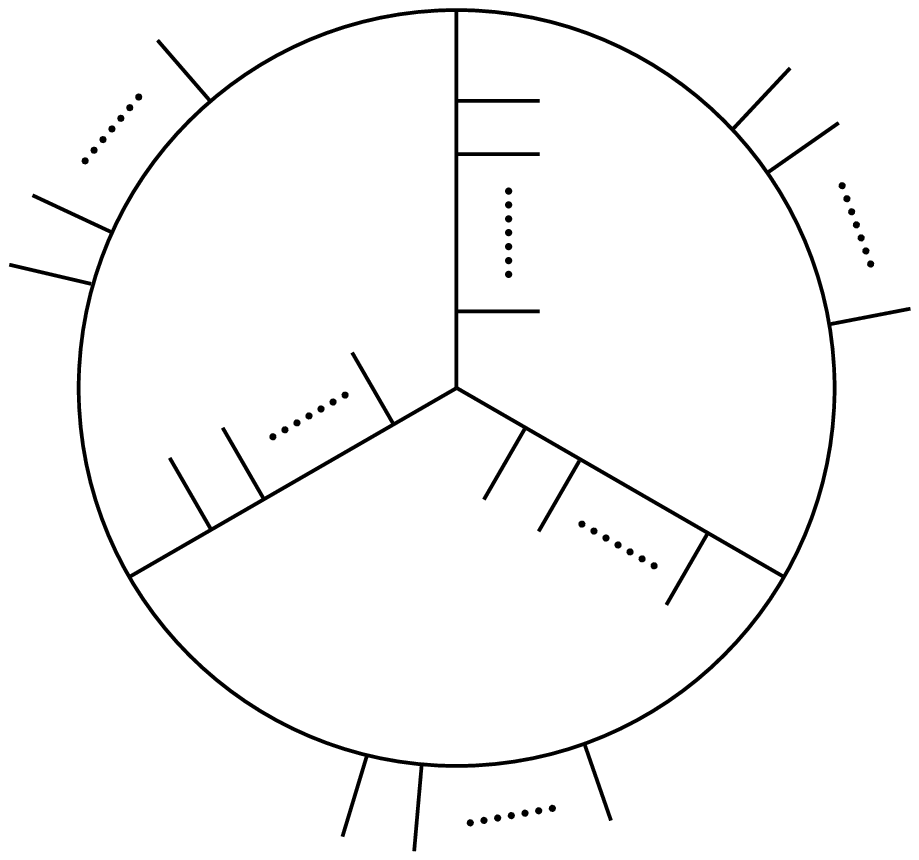}}
\put(-22,70){\footnotesize $n_1$ \!\!\! legs}
\put(88,70){\footnotesize $n_2$ \!\!\! legs}
\put(65,5){\footnotesize $n_3$ \!\!\! legs}
\put(55,60){\footnotesize $n_6$}
\put(55,52){\footnotesize legs}
\put(32,24){\footnotesize $n_4$ \!\!\! legs}
\put(14,53){\footnotesize $n_5$ \!\!\! legs}
\end{picture}\end{array}
\mbox{is identified with }
x_1^{n_1} x_2^{n_2} x_3^{n_3} x_4^{n_4} x_5^{n_5} x_6^{n_6} ,
\end{equation}
and the following relations (as algebra relations) imply the IHX
relations on the legs:
\begin{equation}
\label{eq.Atet_rel}
\begin{cases}
\ x_1 - x_2 - x_6 = 0 , \\
\ x_1 - x_3 + x_5 = 0 , \\
\ x_4 + x_5 + x_6 = 0 .
\end{cases}
\end{equation}
In order to better describe the action of ${\mathfrak S}_4$,
following \cite{Nts98}, we make the substitution
$$
\begin{cases}
\ y_{1}=x_{1}-x_{5}+x_{6} , \\
\ y_{2}=x_{2}+x_{4}-x_{6} , \\
\ y_{3}=x_{3}-x_{4}+x_{5} , \\
\ y_{4}=-x_{1}-x_{2}-x_{3} ,
\end{cases}
$$
replacing variables corresponding with edges of the tetrahedron with
variables corresponding with its faces. In these new variables,
\begin{align*}
\A \big(\tetrahedron \big) & \ \cong \ \Q[y_{1},y_{2},y_{3},y_{4}]
\big/ (y_{1}\!+\!y_{2}\!+\!y_{3}\!+\!y_{4}=0),\mathfrak{S}_{4} \\* &
\ \cong \ \Q[y_{1},y_{2},y_{3},y_{4}]^{{\mathfrak S}_4} \big/
(y_{1}\!+\!y_{2}\!+\!y_{3}\!+\!y_{4}=0),
\end{align*}
where ${\mathfrak S}_4$ acts on $\Q[y_{1},y_{2},y_{3},y_{4}]$ by
permuting $y_1,y_2,y_3,y_4$ symmetrically in even degrees and
skew-symmetrically in odd degrees.

We may identify $\A \big(\tsq\big)$ with the polynomial algebra on
six letters modulo the IHX relations on the legs and modulo the
action of the automorphism group of the $\tsq$--shape as above.
Thus:
$$
\A \big(\tsq\big) \ \cong \ \Q[z_1,z_2,z_3,z_4] \big/
(z_1+z_2+z_3+z_4=0), \mathrm{Aut}(\tsq) ,
$$
where
\begin{equation}
\label{eq.tsq_polyn}
\begin{array}{c}
\begin{picture}(90,90)
\put(0,0){\includegraphics[width=80pt]{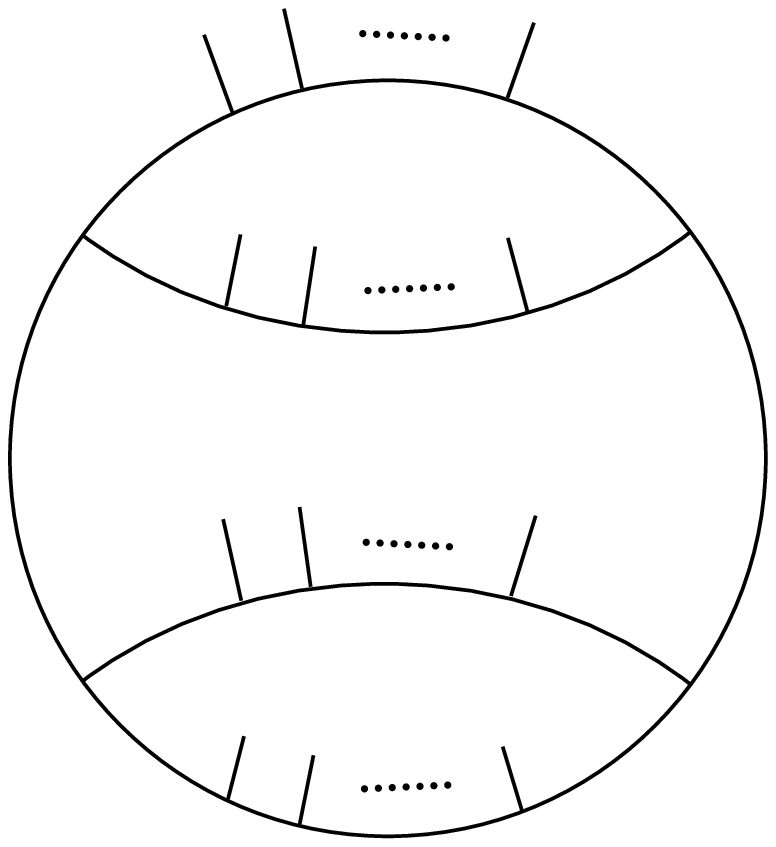}}
\put(26,90){\footnotesize $m_1$ \!\!\! legs}
\put(26,66){\footnotesize $m_2$ \!\!\! legs}
\put(26,37){\footnotesize $m_3$ \!\!\! legs}
\put(26,13){\footnotesize $m_4$ \!\!\! legs}
\end{picture}\end{array}
\mbox{is identified with }
z_1^{m_1} z_2^{m_2} z_3^{m_3} z_4^{m_4} .
\end{equation}

Jacobi diagrams whose internal graphs are $\tsq$ and $\tetrahedron$
are related by the IHX relation
which is obtained from the fourth graph of (\ref{eq.4vg_ihx}),
\begin{equation}
\label{eq.ihx_tsq_tet}
\begin{array}{c}
\begin{picture}(82,95)
\put(0,0){\includegraphics[width=80pt]{tsq1}}
\put(26,90){\footnotesize $m_1$ \!\!\! legs}
\put(26,66){\footnotesize $m_2$ \!\!\! legs}
\put(26,37){\footnotesize $m_3$ \!\!\! legs}
\put(26,13){\footnotesize $m_4$ \!\!\! legs}
\end{picture}\end{array}
\underset{IHX}{=} \
\begin{array}{c}
\begin{picture}(82,95)
\put(0,0){\includegraphics[width=80pt]{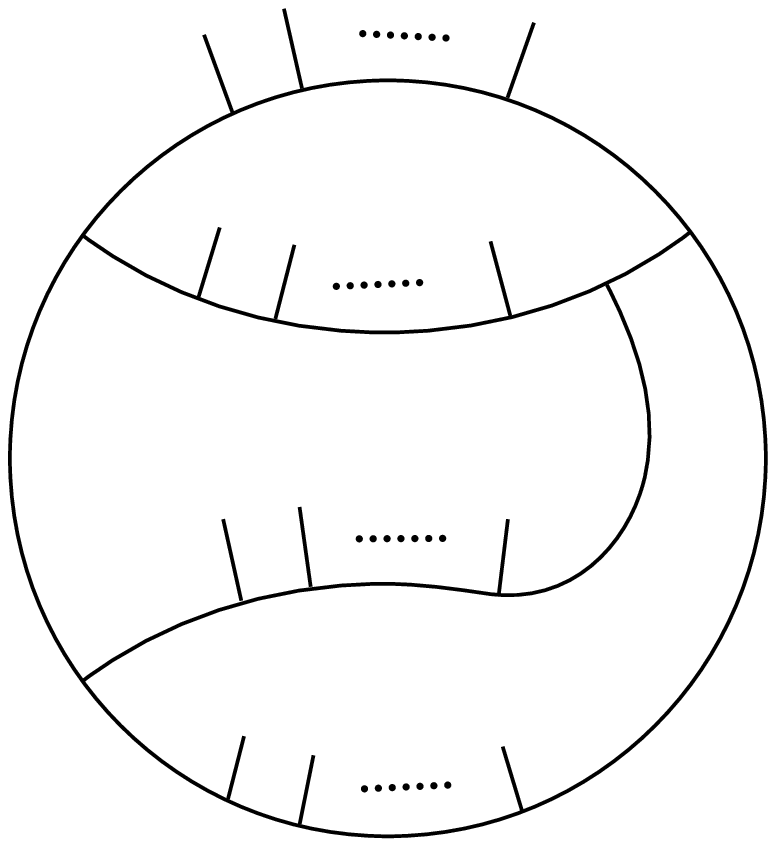}}
\put(26,90){\footnotesize $m_1$ \!\!\! legs}
\put(26,66){\footnotesize $m_2$ \!\!\! legs}
\put(26,37){\footnotesize $m_3$ \!\!\! legs}
\put(26,13){\footnotesize $m_4$ \!\!\! legs}
\end{picture}\end{array}
+ \
\begin{array}{c}
\begin{picture}(85,95)
\put(0,0){\includegraphics[width=80pt]{tet4}}
\put(26,90){\footnotesize $m_2$ \!\!\! legs}
\put(26,66){\footnotesize $m_1$ \!\!\! legs}
\put(26,37){\footnotesize $m_3$ \!\!\! legs}
\put(26,13){\footnotesize $m_4$ \!\!\! legs}
\end{picture}\end{array} .
\end{equation}

\subsection{Odd degree part}
\label{sec.odd_deg}

The aim of this section is to prove the following theorem.

\begin{thm}
\label{thm.main}
The space of $3$--loop Jacobi diagrams of odd degree vanishes.
That is, $\A_{\, \mbox{\scriptsize 3--loop}}^{\rm (odd)} = 0$.
\end{thm}

\begin{proof}
By (\ref{eq.A3l}),
$$
\A_{\, \mbox{\scriptsize 3--loop}}^{\rm (odd)}  \ \cong \
\Big( \bigoplus_{\mbox{\scriptsize $\Gamma$ in (\ref{eq.5ig})}} \!\!
\A (\Gamma)^{\rm (odd)} \Big) \Big/ \, {\rm IHX} .
$$
We show the vanishing of $\A (\Gamma)^{\rm (odd)}$
for the first four graphs $\Gamma$ in (\ref{eq.5ig}).

The vanishing of $\A \big( \tsq \big)^{\rm (odd)}$ is shown as follows.
It is shown by the IHX
relation that this space is spanned by diagrams of the form
(\ref{eq.tsq_polyn}). Such a diagram $D$ is equal modulo the AS
relation to $-D$ by reflection of the internal graph with respect to
a vertical line, therefore $D=0$.
Hence, $\A \big( \tsq \big)^{\rm (odd)} = 0$.

Similarly, reflection of the internal graph shows us that the spaces
$\A \big( \bbl \big)^{\rm (odd)}$ and $\A \big( \mdl \big)^{\rm
(odd)}$ also both vanish.

The vanishing of $\A (\wtr)^{\rm (odd)}\vspc{15}$ is shown as
follows. Let $D$ be a Jacobi diagram whose internal graph is $\wtr$.
We can assume by the IHX relation that there are no legs adjacent to
any separating arc. If there is a loop with an even number of legs,
then the AS relation on the vertex connecting a separating arc with
this loop gives $D=-D$ and therefore $D=0$. Otherwise, by applying
the IHX relation to a separating arc, $D$ is equal to $2$ times a
Jacobi diagram in $\A(\bbl)^{\rm (odd)} = 0$, and therefore $D=0$.
Hence, $\A (\wtr)^{\rm (odd)} =0$.

Therefore, the space of $3$--loop Jacobi diagrams of odd degree is presented by
$$
\A_{\, \mbox{\scriptsize 3--loop}}^{\rm (odd)}  \ \cong \ \A \big(
\tetrahedron \big)^{\rm (odd)} \big/ \big( \mbox{(the right hand
side of (\ref{eq.ihx_tsq_tet}))} = 0 \big).
$$
The vector space spanned by the right hand side of (\ref{eq.ihx_tsq_tet})
is spanned by
$$
\big( x_1^{m_1} x_5^{m_2} + x_1^{m_2} x_5^{m_1} \big) x_4^{m_3} (-x_2)^{m_4}
$$
in terms of polynomials under the identification (\ref{eq.tet_polyn}).
This space is spanned by
$$
(x_1+x_5)^m (x_1 x_5)^n x_4^{m_3} (-x_2)^{m_4} .
$$
Noting that $x_1+x_5 = x_3 = x_2-x_4$,
this space is further spanned by diagrams of the following form.
\begin{equation}
\label{eq.image}
\begin{array}{c}
\begin{picture}(100,90)
\put(0,0){\includegraphics[width=90pt]{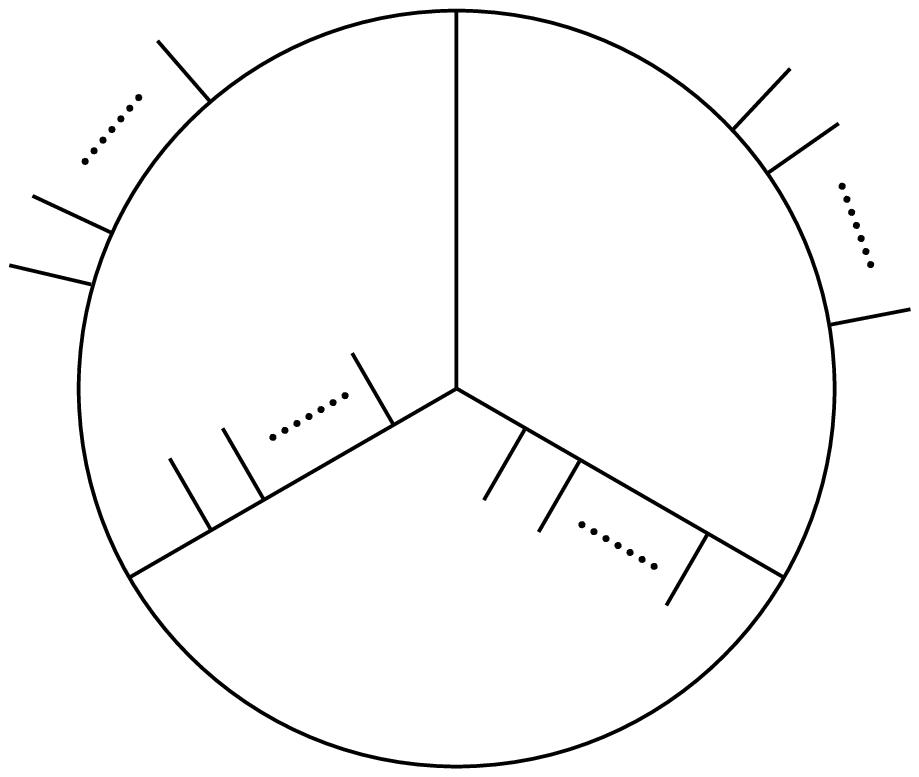}}
\put(-17,65){\footnotesize $n$ \!\!\! legs}
\put(88,65){\footnotesize $n$ \!\!\! legs}
\put(35,15){\footnotesize $n_4$ \!\!\! legs}
\put(15,47){\footnotesize $n_5$ \!\!\! legs}
\end{picture}\end{array}
\end{equation}
Hence,
$$
\A_{\, \mbox{\scriptsize 3--loop}}^{\rm (odd)}  \ \cong \ \A \big(
\tetrahedron \big)^{\rm (odd)} \big/ \big( \mbox{(\ref{eq.image})} =
0 \big) .
$$

In order to show that $\A_{\, \mbox{\scriptsize 3--loop}}^{\rm
(odd)} = 0$, it is sufficient to show that $\A \big( \tetrahedron
\big)^{\rm (odd)}$ is spanned by diagrams of the form
(\ref{eq.image}). As mentioned in Section \ref{sec.polyn_Jd}, the
space of $3$--loop Jacobi diagrams of odd degree is presented by
$$
\A \big(\tetrahedron \big)^{\rm (odd)} \ \cong \ \big(
\Q[y_1,y_2,y_3,y_4]^{\rm (odd)} \big)^{{\mathfrak S}_4} \big/
(y_1\!+\!y_2\!+\!y_3\!+\!y_4=0),
$$
where the action of $\mathfrak{S}_{4}$ on
$\Q[y_{1},y_{2},y_{3},y_{4}]^{\text{(odd)}}$
is skew symmetric.
Since a skew symmetric polynomial is presented by
the product of a symmetric polynomial and
the discriminant $\Delta = \prod_{i<j}(y_{i}-y_{j})$,
$$
\A \big(\tetrahedron \big)^{\rm (odd)}
\ \cong \
\Delta \cdot \Q[\sigma_{2},\sigma_{3},\sigma_{4}]^{\rm (odd)}
\ \cong \
\Delta\sigma_{3}\cdot\Q[\sigma_{2},\sigma_{3}^{2},\sigma_{4}] ,
$$
recalling that
$\sigma_i$ denotes the $i$th symmetric polynomial in $y_1,y_2,y_3,y_4$.
Hence, the vector space spanned by the diagrams of the form (\ref{eq.image})
in $\A \big(\tetrahedron \big)^{\rm (odd)}$
is presented by the image of the following map
$$
\Q[x_1 x_2, x_4, x_5]^{\rm (odd)} \longrightarrow
\Q[x_1,x_2,x_3,x_4,x_5,x_6]^{\rm (odd)} \big/ (\ref{eq.Atet_rel}),
{\mathfrak S}_4 \ \cong \
\Delta\sigma_{3}\cdot\Q[\sigma_{2},\sigma_{3}^{2},\sigma_{4}] .
$$
By Lemma \ref{lem.surj}, this map is surjective, noting that
\begin{align*}
x_1 &= (y_1-y_4)/4 , \\
x_2 &= (y_2-y_4)/4 , \\
x_4 &= (y_2-y_3)/4 , \\
x_5 &= (y_3-y_1)/4 .
\end{align*}
Therefore,
$\A \big( \tetrahedron \big)^{\rm (odd)}$ is spanned by
the diagrams of the form (\ref{eq.image}),
which implies the theorem.
\end{proof}

\begin{rem}
\label{rem.Dasbach_gap} Dasbach \cite{Das97} claimed to have proved
the vanishing of $n$--loop Jacobi diagrams with odd number of legs
for $n\leq 6$ (cited in two of his subsequent papers--- in
\cite{Das98} as Theorem 2.2 and half of Theorem 7.4, and in
\cite{Das00}, although the focus of both papers is on the `even
number of legs' case). There is however a gap in the proof of his
Theorem 5.4.3(\textrm{iii}) (the second equation on page 58 is
wrong, since he is using `modulo greater CW--vectors' to go one way
but not the other).
\end{rem}

\vskip 0.5pc

\begin{rem}
\label{rem.q_inv} It is known that no quantum invariant can
distinguish a knot and its inverse. Hence, if there existed a
counter-example to the conjecture that Jacobi diagrams with an odd
number of legs vanish, such a Jacobi diagram would not be detectable
by weight systems derived from Lie algebras. It is known
\cite{Vogel,Lieberum} how to construct elements which can not be
detected by weight systems derived from Lie algebras, but the method
employed in these papers would not give non-trivial diagrams with an
odd number of legs, as it involves constructing non-trivial diagrams
by multiplying particular elements of Vogel's algebra $\Lambda$, and
the action of $\Lambda$ does not change the number of legs.
\end{rem}

\subsection{Even degree part}
\label{sec.even_deg}

In this section, we review the identification of the space of
$3$--loop Jacobi diagrams of even degree with a polynomial algebra,
following Nakatsuru \cite{Nts98}.
This identification recovers \cite[Theorem 7.4]{Das98}.

By (\ref{eq.A3l}),
$$
\A_{\, \mbox{\scriptsize 3--loop}}^{\rm (even)}  \ \cong \
\Big( \bigoplus_{\mbox{\scriptsize $\Gamma$ in (\ref{eq.5ig})}} \!\!
\A (\Gamma)^{\rm (even)} \Big) \Big/ \, {\rm IHX} .
$$
Unlike the odd degree case, it is necessary to describe IHX
relations among internal graphs $\Gamma$ concretely, since $\A
(\Gamma)^{\rm (even)}$ do not vanish for most $\Gamma$. Let
$\cD(\Gamma)$ denote the space of Jacobi diagrams whose internal
graph is $\Gamma$, not divided by the action of the symmetry of
$\Gamma$. Then, by definition, $\A(\Gamma) = \cD(\Gamma) / {\rm
Aut}\, (\Gamma)$. The IHX relations obtained from the first 4 graphs
of (\ref{eq.4vg_ihx}) induce the maps
\begin{align*}
& \psi_1 : \cD( \wtr ) \longrightarrow \cD( \bbl ), \\
& \psi_2 : \cD( \bbl ) \longrightarrow \cD( \mdl ), \\
& \psi_3 : \cD( \mdl ) \longrightarrow \cD( \tsq ), \\
& \psi_4 : \cD( \tsq ) \longrightarrow \cD( \tetrahedron ).
\end{align*}
Here, for example, $\psi_4$ is the map taking the left hand side of
(\ref{eq.ihx_tsq_tet}) to the right hand side of
(\ref{eq.ihx_tsq_tet}).
Further, the IHX relation obtained from the last graph of (\ref{eq.4vg_ihx})
is the relations,
\begin{equation}
\label{eq.ihx_theta-o}
\begin{picture}(80,30)
\put(0,-2){\pc{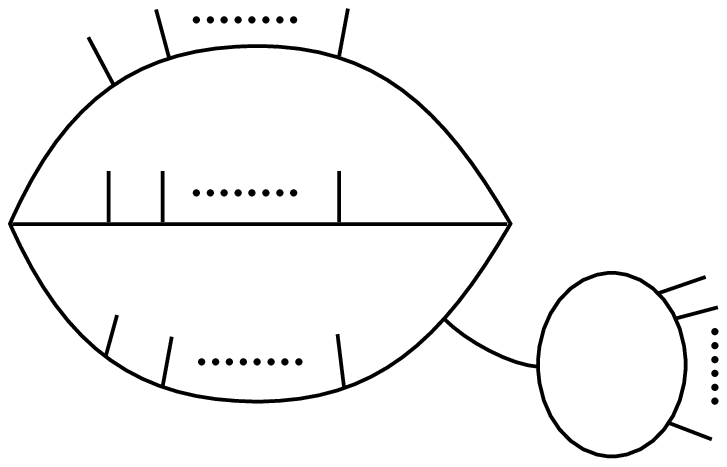}{0.33}}
\put(12,28){\scriptsize $n_1$ legs}
\put(15,11){\scriptsize $n_2$ legs}
\put(17,-4){\scriptsize $n_3$ legs}
\put(75,-14){\scriptsize $n$ legs}
\end{picture}
+
\begin{picture}(80,30)
\put(0,-2){\pc{g10}{0.33}}
\put(12,28){\scriptsize $n_2$ legs}
\put(15,11){\scriptsize $n_3$ legs}
\put(17,-4){\scriptsize $n_1$ legs}
\put(75,-14){\scriptsize $n$ legs}
\end{picture}
+
\begin{picture}(80,30)
\put(0,-2){\pc{g10}{0.33}}
\put(12,28){\scriptsize $n_3$ legs}
\put(15,11){\scriptsize $n_1$ legs}
\put(17,-4){\scriptsize $n_2$ legs}
\put(75,-14){\scriptsize $n$ legs}
\end{picture}
= \ 0 .
\end{equation}
By using these,\vspc{23}
the space of $3$--loop Jacobi diagrams of even degree is presented by
$$
\A_{\, \mbox{\scriptsize 3--loop}}^{\rm (even)}  \ \cong \
\Big( \bigoplus_{\mbox{\scriptsize $\Gamma$ in (\ref{eq.5ig})}} \!\!
\cD (\Gamma)^{\rm (even)} \Big) \Big/ \,
\big( {\rm Aut}\,(\Gamma) \mbox{ for $\Gamma$ in (\ref{eq.5ig})}, \
\psi_1,\, \psi_2, \, \psi_3, \, \psi_4, \, \mbox{(\ref{eq.ihx_theta-o})} \big).
$$

Since $\A( \wtr )^{\rm (even)} = 0$ and $\psi_1$ induces the zero
map $\A( \wtr )^{\rm (even)} \to \A( \mdl )^{\rm (even)}$, we can
ignore the contribution from $\A( \wtr )^{\rm (even)}$.
Further, since $\psi_3 \psi_2$ descends to
a map $\A( \bbl )^{\rm (even)} \to \A( \tsq )^{\rm (even)}$,
we can ignore the contribution from $\A( \bbl )^{\rm (even)}$.
Furthermore, since $\psi_3$ induces
a map $\A( \mdl )^{\rm (even)} \to \A( \tsq )^{\rm (even)}$ and
(\ref{eq.ihx_theta-o}) vanishes in the image of $\psi_4 \psi_3$, we
can ignore the contribution from $\A( \mdl )^{\rm (even)}$. Hence,
$$
\A_{\, \mbox{\scriptsize 3--loop}}^{\rm (even)}  \ \cong \
\Big( \cD (\tsq)^{\rm (even)} \oplus \cD (\tetrahedron)^{\rm (even)} \Big)
\Big/
\Big( {\rm Aut}\,(\tsq), \, {\rm Aut}\,(\tetrahedron), \, \psi_4 \Big) .
$$

It can be checked by concrete calculation that
if Jacobi diagrams $D, D' \in \cD \big(\tsq \big)^{\rm (even)}$
are related by ${\rm Aut}\, \big( \tsq \big)$,
then $\psi_4(D)$ and $\psi_4(D')$ are related by
${\rm Aut}\, \big( \tetrahedron \big)$.
Hence, $\psi_4$ induces
a map $\overline{\psi_4} :
\A( \tsq )^{\rm (even)} \to \A( \tetrahedron )^{\rm (even)}$.
Therefore,
$$
\A_{\, \mbox{\scriptsize 3--loop}}^{\rm (even)}  \ \ \cong \ \
\Big( \A (\tsq)^{\rm (even)} \oplus \A (\tetrahedron)^{\rm (even)} \Big)
\Big/ \, \overline{\psi_4}
\ \ \cong \ \
\A (\tetrahedron)^{\rm (even)} .
$$

Hence, by the identification of $\A (\tetrahedron)$ with the
polynomial algebra mentioned in Section \ref{sec.polyn_Jd},
\begin{align*}
\A_{\, \mbox{\scriptsize 3--loop}}^{\rm (even)} & \ \cong \ \big(
\Q[y_{1},y_{2},y_{3},y_{4}]^{\rm (even)}
\big)^{\mathfrak{S}_{4}}\big/
(y_{1}\!+\!y_{2}\!+\!y_{3}\!+\!y_{4}=0) \\
& \ \cong \
\Q[\sigma_{2},\sigma_{3},\sigma_{4}]^{\rm (even)}
\ \cong \ \Q[\sigma_{2},\sigma_{3}^{2},\sigma_{4}] ,
\end{align*}
where $\sigma_{i}$ denotes the $i$th symmetric polynomial
in four variables $y_1,y_2,y_3,y_4$.
It has as its generating function
$$
\frac{1}{(1-x^{2})(1-x^{4})(1-x^{6})}=\sum_{n \text{
even}}\left(\Bigl\lfloor\frac{n^{2}+12n}{48}\Bigr\rfloor+1\right)x^{n}
$$
recovering \cite[Theorem 7.4]{Das98} and agreeing with the results
of \cite{Das00}.

\section{A lemma on polynomial algebras}
\label{sec.pf_lem}

The aim of this section is to prove Lemma \ref{lem.surj}, which was
used in the proof of the main theorem in the previous section.

The {\it skew symmetrizer}
$$
\Q[y_1,y_2,y_3,y_4] \ \longrightarrow \
\Delta \cdot \Q[\sigma_1,\sigma_2,\sigma_3,\sigma_4]
$$
is the linear map sending
$$
f(y_1,y_2,y_3,y_4) \ \mbox{ to } \ \frac{1}{4!} \sum_{\tau \in
{\mathfrak S}_4} \mbox{sgn}(\tau) \,
f(y_{\tau(1)},y_{\tau(2)},y_{\tau(3)},y_{\tau(4)}),
$$
where $\sigma_i$ is the $i$th symmetric polynomial in
$y_1,y_2,y_3,y_4$ and $\Delta = \prod_{i<j}(y_{i}-y_{j})$ as before.
We consider the composition
\begin{multline*} \Q[y_1-y_3, y_2-y_3,
(y_1-y_4)(y_2-y_4)] \ \longrightarrow \\ \Q[y_1,y_2,y_3,y_4]/(y_1 +
y_2 + y_3 + y_4) \ \longrightarrow \ \Delta \cdot
\Q[\sigma_2,\sigma_3,\sigma_4]
\end{multline*}
where the first map is the projection of the inclusion,
and the second map is a quotient of the skew symmetrizer.

\begin{lem}
\label{lem.surj}
The odd degree part of the above map,
$$
\Q[y_1 \!-\! y_3, \, y_2 \!-\! y_3, \,
(y_1 \!-\! y_4)(y_2 \!-\! y_4)]^{\rm (odd)} \ \longrightarrow \
\Delta \sigma_3 \cdot \Q[\sigma_2,\sigma_3^2,\sigma_4],
$$
is surjective,
where $\Q[\cdots]^{\rm (odd)}$ denotes the vector subspace of $\Q[\cdots]$
spanned by polynomials of odd degrees.
\end{lem}

\begin{proof}
We put
\begin{align*}
& P_2(y_1,y_2,y_3) = (y_1-y_2)^2 + (y_2-y_3)^2 + (y_3-y_1)^2 , \\
& P_3(y_1,y_2,y_3) = (y_1-y_2)(y_2-y_3)(y_3-y_1) , \\
& P_4(y_1,y_2,y_3,y_4) = (y_1-y_3)(y_2-y_3)(y_1-y_4)(y_2-y_4) .
\end{align*}
By definition,
$$
12 \, P_2(y_1,y_2,y_3)^n \, P_3(y_1,y_2,y_3)^{2m+3} \, P_4(y_1,y_2,y_3,y_4)^k
$$
belongs to $\Q[y_1\!-\!y_3, \, y_2\!-\!y_3, \,
(y_1\!-\!y_4)(y_2\!-\!y_4)]^{\rm (odd)}$ for any non-negative
integers $n,m,k$. Since $P_2(y_1,y_2,y_3)$ and $P_3(y_1,y_2,y_3)$
are invariant under cyclic permutations of $y_1,y_2,y_3$, the above
polynomial and
\begin{align*}
& 4 \, P_2(y_1,y_2,y_3)^n \, P_3(y_1,y_2,y_3)^{2m+3}  \\*
& \quad \times\big(
P_4(y_1,y_2,y_3,y_4)^k + P_4(y_1,y_3,y_2,y_4)^k + P_4(y_1,y_4,y_2,y_3)^k \big)
\end{align*}
are taken to the same image by the skew symmetrizer.
Further, since the last factor of the above formula is a symmetric polynomial,
the skew symmetrizer takes the above formula to
\begin{align*}
Q^{n,m,k} & = \big( P_2(y_1,y_2,y_3)^n P_3(y_1,y_2,y_3)^{2m+3}
+ P_2(y_4,y_3,y_2)^n P_3(y_4,y_3,y_2)^{2m+3} \\*
& \quad + P_2(y_3,y_4,y_1)^n P_3(y_3,y_4,y_1)^{2m+3}
+ P_2(y_2,y_1,y_4)^n P_3(y_2,y_1,y_4)^{2m+3} \big) \\*
& \times \big(
P_4(y_1,y_2,y_3,y_4)^k + P_4(y_1,y_3,y_2,y_4)^k + P_4(y_1,y_4,y_2,y_3)^k \big).
\end{align*}
Hence, it is sufficient to show that
$\Delta \sigma_3 \cdot \Q[\sigma_2,\sigma_3^2,\sigma_4]$
is spanned by $Q^{n,m,k}$.

For a fixed non-negative integer $d$, we consider the vector
subspace of $\Delta \sigma_3 \cdot \Q[\sigma_2,\sigma_3^2,\sigma_4]$
spanned by polynomials of degree $2d+9$. Since it is spanned by
$\Delta \sigma_3 \cdot \sigma_2^n \sigma_3^{2m} \sigma_4^k$ for
non-negative integers $n,m,k$ satisfying $n+2k+3m=d$, its dimension
is equal to the number of such $(n,m,k)$. Since the $Q^{n,m,k}$'s
are such polynomials of this number, it is sufficient to show the
linear independence of $Q^{n,m,k}$ for non-negative integers $n,m,k$
satisfying that $n+2k+3m=d$.

In order to prove the linear independence of the $Q^{n,m,k}$'s we
first make the substitution
\begin{align*}
y_1 &= (3 t^a - t^b - t^c)/4, \\
y_2 &= (- t^a +3 t^b - t^c)/4, \\
y_3 &= (- t^a - t^b +3 t^c)/4, \\
y_4 &= -(t^a + t^b + t^c)/4,
\end{align*}
where $t$ is a variable tending to $\infty$,
and $a,b,c$ are real numbers satisfying that
$a>b>c>0$ and $a-b \, < \, b-c \, < \, 2(a-b)$.
Since
\begin{align*}
& y_1-y_4 = t^a, \\
& y_2-y_4 = t^b, \\
& y_3-y_4 = t^c,
\end{align*}
we have that
\begin{align*}
P_2(y_1,y_2,y_3) &=
(t^a-t^b)^2 + (t^b-t^c)^2 + (t^c-t^a)^2 \\*
&= 2 t^{2a} \big( 1 - t^{-(a-b)} + o(t^{-(b-c)}) \big),
\end{align*}
where $f(t)=g(t)+o(t^\eps)$ means that
$\big( f(t)-g(t) \big)/t^\eps \to 0$ as $t \to \infty$.
Hence,
$$
P_2(y_1,y_2,y_3)^n = 2^n t^{2a n} \big( 1 - n t^{-(a-b)} + o(t^{-(b-c)}) \big).
$$
Similarly,
\begin{align*}
&P_2(y_4,y_3,y_2)^n = 2^n t^{2b n} \big( 1 + o(t^{0}) \big), \\
&P_2(y_3,y_4,y_1)^n = 2^n t^{2a n} \big( 1 + o(t^{0}) \big), \\
&P_2(y_2,y_1,y_4)^n=2^n t^{2a n}\big(1-n\, t^{-(a-b)} + o(t^{-(b-c)}) \big), \\
&P_3(y_1,y_2,y_3)^{2m+3} = -t^{(2a+b)(2m+3)}
   \big(1- (2m+3) (t^{-(a-b)}+t^{-(b-c)}) + o(t^{-(b-c)}) \big), \\
&P_3(y_4,y_3,y_2)^{2m+3} = t^{(2b+c)(2m+3)} \big( 1 + o(t^{0}) \big), \\
&P_3(y_3,y_4,y_1)^{2m+3} = -t^{(2a+c)(2m+3)} \big( 1 + o(t^{0}) \big), \\
&P_3(y_2,y_1,y_4)^{2m+3} = t^{(2a+b)(2m+3)}
   \big(1- (2m+3) \, t^{-(a-b)} + o(t^{-(b-c)}) \big), \\
&P_4(y_1,y_2,y_3,y_4)^k = t^{(2a+2b)k} \big( 1 + o(t^{0}) \big), \\
&P_4(y_1,y_3,y_2,y_4)^k = (-1)^{k}t^{(2a+b+c)k} \big( 1 + o(t^{0}) \big), \\
&P_4(y_1,y_4,y_2,y_3)^k = t^{(2a+b+c)k} \big( 1 + o(t^{0}) \big).
\end{align*}
Hence,
\begin{align*}
& P_2(y_1,y_2,y_3)^n P_3(y_1,y_2,y_3)^{2m+3}
+ P_2(y_4,y_3,y_2)^n P_3(y_4,y_3,y_2)^{2m+3} \\*
& \quad + P_2(y_3,y_4,y_1)^n P_3(y_3,y_4,y_1)^{2m+3}
+ P_2(y_2,y_1,y_4)^n P_3(y_2,y_1,y_4)^{2m+3} \\
& = -2^n t^{2a n+(2a+b)(2m+3)} \big(1- (n+2m+3) \, t^{-(a-b)}
-(2m+3) \, t^{-(b-c)} + o(t^{-(b-c)}) \big) \\* & \quad +2^n t^{2b
n+(2b+c)(2m+3)} \big( 1 + o(t^{0}) \big) \\* & \quad -2^n t^{2a
n+(2a+c)(2m+3)} \big( 1 + o(t^{0}) \big) \\* & \quad +2^n t^{2a
n+(2a+b)(2m+3)} \big(1- (n+2m+3) \, t^{-(a-b)} + o(t^{-(b-c)}) \big)
\\* & = 2^n (2m+3) \, t^{2a n+(2a+b)(2m+3)-(b-c)} \big( 1 + o(t^{0})
\big),
\end{align*}
noting that we need $2m+3>1$ when we verify that
$$
2a n+(2a+b)(2m+3)-(b-c) \ > \ 2a n+(2a+c)(2m+3).
$$
Further,
$$
P_4(y_1,y_2,y_3,y_4)^k + P_4(y_1,y_3,y_2,y_4)^k + P_4(y_1,y_4,y_2,y_3)^k
= \eps \, t^{(2a+2b)k} \big( 1 + o(t^{0}) \big) ,
$$
where $\eps=1$ if $k>0$, and $\eps=3$ if $k=0$.
Therefore,
$$
Q^{n,m,k} = \eps \, 2^n (2m+3) \, t^{ 2(n+2m+k+3)a + 2(m+k+1)b +c }
\big( 1 + o(t^{0}) \big) .
$$
This implies that the only possible linear relations between the
$Q^{n,m,k}$'s are between those having the same value of $(n+2m+k, \
m+k )$, and in particular the same value of $m+k$. In other words,
the vector space which we are considering is presented by the direct
sum:
$$
\mbox{span} \{ Q^{n,m,k} \ | \ n+2k+3m=d \} \ = \
\bigoplus_{\ell}
\mbox{span} \{ Q^{n,m,k} \ | \ n+2k+3m=d, \ \ m+k=\ell \} .
$$

Next, in order to complete the proof of the linear independence of the
$Q^{n,m,k}$'s,
we make another substitution
\begin{align*}
y_1 &= (2 t^a - t^c)/4\thickspace\ + t^b/2 , \\
y_2 &= (2 t^a - t^c)/4\thickspace\  - t^b/2 , \\
y_3 &= (-2 t^a + 3 t^c)/4 , \\
y_4 &= -(2 t^a + t^c)/4 ,
\end{align*}
where $t$ is as above,
and $a,b,c$ are real numbers satisfying that
$a>b>c>0$ and $b-c \, < \, a-b \, < \, 2(b-c)$.
Since
\begin{align*}
& y_1-y_4 = t^a + t^b/2, \\
& y_2-y_4 = t^a - t^b/2, \\
& y_3-y_4 = t^c , \\
& y_1-y_2 = t^b ,
\end{align*}
we have that
\begin{align*}
&P_2(y_1,y_2,y_3)^n = 2^n t^{2a n} \big( 1 -2 n \, t^{-(a-c)} + o(t^{-(a-c)}) \big), \\
&P_2(y_4,y_3,y_2)^n = 2^n t^{2a n}
     \big( 1 -n\, t^{-(a-b)} -n\, t^{-(a-c)}+ o(t^{-(a-c)}) \big), \\
&P_2(y_3,y_4,y_1)^n = 2^n t^{2a n}
     \big( 1 +n\, t^{-(a-b)} -n\, t^{-(a-c)} + o(t^{-(a-c)}) \big), \\
&P_2(y_2,y_1,y_4)^n = 2^n t^{2a n} \big( 1 + o(t^{-(a-c)}) \big), \\
&P_3(y_1,y_2,y_3)^{2m+3} = -t^{(2a+b)(2m+3)}
   \big(1- 2 (2m+3) \, t^{-(a-c)} + o(t^{-(a-c)}) \big), \\
&P_3(y_4,y_3,y_2)^{2m+3} = t^{(2a+c)(2m+3)}
     \big( 1 - (2m+3)\, (t^{-(a-b)} +  t^{-(a-c)}) + o(t^{-(a-c)}) \big), \\
&P_3(y_3,y_4,y_1)^{2m+3} = -t^{(2a+c)(2m+3)}
     \big( 1 + (2m+3)\, (t^{-(a-b)} - t^{-(a-c)}) + o(t^{-(a-c)}) \big), \\
&P_3(y_2,y_1,y_4)^{2m+3} = t^{(2a+b)(2m+3)} \big( 1 + o(t^{-(a-c)}) \big), \\
&P_4(y_1,y_2,y_3,y_4)^k = t^{4a k} \big( 1 + o(t^{0}) \big), \\
&P_4(y_1,y_3,y_2,y_4)^k = (-1)^{k}t^{(2a+b+c)k} \big( 1 + o(t^{0}) \big), \\
&P_4(y_1,y_4,y_2,y_3)^k = t^{(2a+b+c)k} \big( 1 + o(t^{0}) \big).
\end{align*}
Hence,
\begin{align*}
& P_2(y_1,y_2,y_3)^n P_3(y_1,y_2,y_3)^{2m+3}
+ P_2(y_4,y_3,y_2)^n P_3(y_4,y_3,y_2)^{2m+3} \\*
& \quad + P_2(y_3,y_4,y_1)^n P_3(y_3,y_4,y_1)^{2m+3}
+ P_2(y_2,y_1,y_4)^n P_3(y_2,y_1,y_4)^{2m+3} \\
& = -2^n t^{2a n+(2a+b)(2m+3)} \big(1 - 2(n+2m+3) \, t^{-(a-c)} +
o(t^{-(a-c)}) \big) \\* & \quad +2^n t^{2a n+(2a+c)(2m+3)} \big(1 -
(n+2m+3) \, (t^{-(a-b)} +  t^{-(a-c)}) + o(t^{-(a-c)}) \big)
\\* & \quad -2^n t^{2a n+(2a+c)(2m+3)} \big(1 + (n+2m+3) \,
(t^{-(a-b)} -  t^{-(a-c)}) + o(t^{-(a-c)}) \big) \\* & \quad +2^n
t^{2a n+(2a+b)(2m+3)} \big(1 + o(t^{-(a-c)}) \big) \\* & = 2^{n+1}
(n+2m+3) \, t^{2a n+(2a+b)(2m+3)-(a-c)} \big( 1 + o(t^{0}) \big).
\end{align*}
Further,
$$
P_4(y_1,y_2,y_3,y_4)^k + P_4(y_1,y_3,y_2,y_4)^k + P_4(y_1,y_4,y_2,y_3)^k
= \eps \, t^{4a k} \big( 1 + o(t^{0}) \big) ,
$$
where $\eps=1$ if $k>0$, and $\eps=3$ if $k=0$.
Therefore,
$$
Q^{n,m,k} =
\eps \, 2^{n+1} (n+2m+3) \, t^{(2(n+2m+2k)+5)a +(2m+3) b +c}
\big( 1 + o(t^{0}) \big).
$$
This implies that the only possible linear relations between the
$Q^{n,m,k}$'s are between those having the same value of $(n+2m+2k, \ m )$.

Thus, the only linear relations
that could exist between $Q^{n,m,k}$'s with fixed $n+2k+3m=d$ are
between those having the same value $m+k$ (from the first substitution)
and the same values of $m$ (from the second substitution).
In other words,
$$
\mbox{span} \{ Q^{n,m,k} \ | \ n+2k+3m=d \} \ = \
\!\!\! \bigoplus_{n+2k+3m=d} \!\!\!
\mbox{span} \{ Q^{n,m,k} \} .
$$
It follows that the $Q^{n,m,k}$'s are indeed linearly independent,
as required.
\end{proof}

\bibliographystyle{amsplain}

\providecommand{\bysame}{\leavevmode\hbox to3em{\hrulefill}\thinspace}

\end{document}